\documentclass[11pt,fleqn,titlepage]{article}
\usepackage{latexsym,amsmath,amssymb,amsfonts,epsfig}
\usepackage[all]{xy}

\providecommand{\RR}{{\mathbb{R}}}

\providecommand{\ZZ}{{\mathbb{Z}}}

\providecommand{\QED}{{\hfill{} $\Box$}}

\providecommand{\GG}{{\cal G}}
\providecommand{\HH}{{\cal H}}
\providecommand{\KK}{{\cal K}}
\providecommand{\PP}{{\mathbb P}}

\providecommand{\THM}{{\mathbb{T}_H M}}
\providecommand{\TM}{{\mathbb{T}M}}

\providecommand{\PP}{{\mathbb{P}}}

\providecommand{\sm}{{C^\infty}}
\newcommand{\ang}[1]{\langle #1 \rangle} 

\newtheorem{definition}{Definition}

\newtheorem{proposition}[definition]{Proposition}
\newtheorem{corollary}[definition]{Corollary}
\newtheorem{theorem}[definition]{Theorem}

\begin{document}

\title{The Index of Hypoelliptic Operators on Foliated Manifolds}
\author{Erik van Erp, Dartmouth College, erik.van.erp@dartmouth.edu}
\maketitle

{\small 
\begin{center}{\bf Abstract}\end{center}
In [vE1] and [vE2] we presented the solution to the index problem
for a class of hypoelliptic operators on closed contact manifolds.
The proofs are based on an adaptation of the tangent groupoid method of Alain Connes to hypoelliptic index problems.
The methods originally developed for contact manifolds have wider applicability to the index theory of hypoelliptic Fredholm operators.
As an illustration of the scope and effectiveness of these methods, we present here an index theorem for a class of hypoelliptic differential operators on closed foliated manifolds.

}

\section{Introduction}


In [vE1], [vE2] we published the solution of the index problem for a class of hypoelliptic (pseudo)differential operators on contact manifolds. 
Our solution relied on techniques from noncommutative geometry,
and makes liberal use of $C^*$-algebras, groupoids, and analytic $K$-theory.
The proof of our index theorem proceeded in two stages. 
The key insight in [vE1] is that the {\em principal part} of a hypoelliptic operator (suitably interpreted in the sense of the {\em Heisenberg calculus}) gives rise to a class in the analytic $K$-theory of a noncommutative algebra.
As it turns out, this analytic $K$-theory group is canonically isomorphic to the topological group $K^0(T^*M)$,
and the main theorem in [vE1] states that the Fredholm index of this class of hypoelliptic operators can be computed by the Atiyah-Singer index formula.

This is, however, not yet a useful theorem, because a class in the analytic $K$-theory of a noncommutative algebra is very hard to compute in concrete examples.
To obtain a cohomological index formula 
we must find a {\em topological} expression for this $K$-theory class as an element in $K^0(T^*M)$.
In [vE2] we solve this problem for contact manifolds and derive a concrete topological construction of a symbol class for hypoelliptic operators in $K^0(T^*M)$.

While both papers [vE1] [vE2]  discuss hypoelliptic index theory for contact manifolds, the results presented in the first of these two papers [vE1] can be stated in much greater generality.
Contact manifolds can be thought of as a specific type of {\em filtered manifolds}, i.e., manifolds equipped with a distribution $H\subseteq TM$ (not necessarily of codimension 1).
A typical hypoelliptic operator of the class we are interested in would be a second order ``elliptic'' operator in the directions of $H$,
but would be of first order transversally to $H$.
More precisely, we study operators that are ``elliptic'' in the Heisenberg calculus associated to $H$.
The ideas and techniques developed in [vE1] can be generalized to apply to all such hypoelliptic operators, regardless of the geometric nature of the structure $H$.
One could thus derive a generalized index theorem for hypoelliptic operators and investigate what it says for other types of geometric structure $H$.

The result of this new approach to hypoelliptic index theory is two-fold.
On the one hand we have a very general theorem (proven along the lines of [vE1]) whose essential import is that the formula of Atiyah and Singer computes the Fredholm index not only for elliptic operators, but for all hypoelliptic operators in the Heisenberg calculus for some bundle $H$.
The index theorem of Boutet de Monvel for Toeplitz operators [Bo] is an example and a special case of this fact: the Atiyah-Singer formula applies to Fredholm Toeplitz operators.
The situation is analogous (but different) from an early theorem by H\"ormander,  who extended the Atiyah-Singer formula to hypoelliptic operators of type $(\rho, \delta)$, with $0\le 1-\rho\le \delta<\rho\le 1$ (See [H\"o]).
The operators covered by our methods are of type $(\frac{1}{2}, \frac{1}{2})$ and therefore not covered by H\"ormander's theorem.
But in both cases the Atiyah-Singer formula computes the index.

The second aspect of our approach  is that in order to apply the topological formula of Atiyah-Singer  to non-elliptic operators we will have to modify the definition of the $K$-theory class associated to the symbol of the operator.
This is the general lesson to be learned from the results of the second paper [vE2]:
to turn the `non-commutative' version of our hypoelliptic index theorem into an explicit cohomological formula one must solve a nontrivial problem, namely to find a topological expression for an analytically constructed $K$-theory class.
In [vE2] we solved this problem if $H$ is a contact structure,
but the ideas used in [vE2] do not apply in the general case.
In fact, there is no general method that solves this problem once and for all for all possible structures $H\subseteq TM$.
(In a forthcoming publication we will discuss this general problem in more detail.)
The details of the calculation of the correct symbol class in $K^0(T^*M)$
depend on specific geometric properties of the structure $H$,
and creative ideas are needed to solve it in each particular case. 

The present paper contains the solution of the hypoelliptic index problem in the Heisenberg calculus in the case where $H\subset TM$ is a foliation.
As it turns out, this problem is easier to handle than the corresponding problem for contact manifolds.
The exposition in this paper is relatively self-contained, and we derive our index theorem, as well as the necessary hypoelliptic Fredholm theory, from scratch.

\section{The Index Formula for Foliations}

Throughout this paper $M$ denotes a smooth closed manifold of dimension $n$.
Consider a differential operator $P$ on $M$ acting on smooth sections in a vector bundle $E$, with values in a second bundle $F$,
\[ P\;\colon\; \sm(E)\to \sm(F).\]
Effective calcultations with such an operator are performed in a coordinate system 
\[ x=(x_1, \ldots, x_n)\;\colon\; U\to \RR^n\]
in an open set $U\subset M$,  where the bundles $E, F$ are trivialized as $U\times \RR^m$ (we will assume in what follows that $E$ and $F$ have equal rank $m$).
Then $P$ is represented as
\[ P = \sum_{|\alpha|\le d} a_\alpha(x) \partial^\alpha.\]
As is customary, $\alpha=(\alpha_1, \ldots ,\alpha_n)$ denotes a multi-index of nonnegative integers;
the coefficients $a_\alpha(x)$ are matrix-valued smooth functions on $U$;
and 
\[ \partial^\alpha = \left(\frac{\partial}{\partial x_1}\right)^{\alpha_1} \cdots 
\left(\frac{\partial}{\partial x_n}\right)^{\alpha_n} .\]
Finally, the degree of the monomial $a_\alpha \partial^\alpha$ is denoted as
\[ |\alpha| = \alpha_1+\cdots +\alpha_n,\]
and the order $d$ of $P$ is the maximal degree $|\alpha|$ for which $a_\alpha\ne 0$ at some point of $M$.

The principal symbol of the operator $P$ is the function
\[ \sigma(P)(x,\xi) = \sum_{|\alpha|=d} a_\alpha(x) (i\xi)^\alpha,\]
for $(x,\xi)\in \RR^n\times \RR^n$.
The principal symbol transforms (under a change of coordinates) as a smooth section in the bundle of algebras
\[ {\rm Hom}(\pi^*E,\pi^*F),\]
over the cotangent space $T^*M$, where $\pi\colon T^*M\to M$ denotes the base point map.
Obviously, $\sigma(P)$ is a homogeneous polynomial in the fibers of $T^*M$. 

The operator $P$ is called {\em elliptic} if the principal symbol is invertible at all points of $M$ (except, of course, at $\xi=0$), and elliptic operators are Fredholm if $M$ is closed.
We will treat $P$ as an unbounded Hilbert space operator 
\[ P\;\colon\; L^2(E)\to L^2(F),\]
with domain $C^\infty(E)$.
The differential operator has a closure, which we denote by $\bar{P}$.
In this context, to say that $P$ is Fredholm means that the closure $\bar{P}$ has closed range and finite dimensional kernel and cokernel. The Fredholm index is defined as
\[ {\rm Index}\, P = {\rm dim Ker}\,P -{\rm dim Coker}\, P,\]
and the formula of Atiyah and Singer [AS1, AS3] computes the index of $P$ as a function of the homotopy type of the principal symbol,
\[ {\rm Index}\,P = \int_{T^*M} {\rm Ch}(\sigma(P))\wedge {\rm Td}(M).\] 
Now let $M$ be foliated by an integrable sub-bundle $H\subset TM$ of the tangent bundle.
We denote the rank of $H$ by $p$, and the rank of the quotient bundle $N=TM/H$ by $q$ (so that $n=p+q$).
When working with a foliation we will always choose coordinates $x$ such that the bundle $H$ (restricted to $U\subset M$) is spanned by the first $p$ coordinate vectorfields $\partial_1, \ldots, \partial_p$.
We now change the usual calculus of differential operators by defining the {\em weighted order} of a monomial $a_\alpha\partial^\alpha$ to be
\[ ||\alpha|| = \alpha_1+\cdots +\alpha_p + 2(\alpha_{p+1}+\cdots +\alpha_n).\]
This amounts to assigning weight two to vector fields that are transversal to the foliation. 
The point of this alternative calculus is to change the notion of {\em highest order part}, or, what amounts to the same thing, the principal symbol of the operator $P$.
With the above notation, the {\em weighted principal symbol} of an operator $P$ of weighted order $d$ is defined in the obvious way as
\[ \sigma_H(P)(x,\xi) = \sum_{||\alpha||= d} a_\alpha(x) (i\xi)^\alpha.\]
It is easily verified that $\sigma_H(P)$ transforms in the same way as $\sigma(P)$, i.e., as a smooth section of the vector bundle
\[ {\rm Hom}(\pi^*E,\pi^*F).\]
The reader who actually checks this last statement will notice that $\sigma_H(P)$ is, canonically, a section in the pull-back bundle over the space $H^*\oplus N^*$ instead of over $T^*M$. Of course, by choosing an arbitary section $N\to TM$ we may identify $T^*M\cong H^*\oplus N^*$, and think if $\sigma_H(P)$ simply as a section over $T^*M$.
While this introduces some arbitrariness, the choice of section $N\to TM$ does not affect the {\em homotopy type} of $\sigma_H(P)$.

Also, for what follows it is worth remarking that $\sigma_H(P)$ is homogeneous of degree $d$ in $\xi$.
Of course, the appropriate notion of homogeneity is the one associated to the grading of the bundle $H\oplus N$, where vectors in $H$ have degree $1$, and vectors in $N$ have degree $2$.

\vskip 12pt
\noindent We can now state our result.

\begin{theorem}\label{thmHormander}
Let $P$ be a differential operator on a closed foliated manifold $(M,H)$.
Suppose the weighted principal symbol $\sigma_H(P)$ is invertible (when $\xi\ne 0$). 
Then $P$ is a Fredholm operator and the Fredholm index of $P$ is computed by the Atiyah-Singer formula,
\[ {\rm Index}(P) = \int_{T^*M} {\rm Ch}(\sigma_H(P))\wedge {\rm Td}(M).\] 
\end{theorem}
The remainder of this paper is devoted to a new proof of this result using the tangent groupoid methodology. The paper is largely self-contained, and we develop the necessary hypoelliptic theory from scratch.

\section{Fredholm Theory for Subelliptic Operators}

One way to prove an operator $P$ is Fredholm is to exhibit an explicit parametrix.
This is the usual approach, and it requires the elaboration of an appropriate pseudodifferential calculus in which the parametrix exists.
In this section we give a more elementary proof of the fact that differential operators with invertible $\sigma_H(P)$ are Fredholm,
without recourse to the full pseudodifferential Heisenberg calculus (see [BG] or [Ta]).  
The basic idea is that one can prove so-called {\em a-priori estimates} for $P$ directly, without first constructing a parametrix.
No matter how it is dressed up, the crucial element in the proof of Fredholmness is the Fourier transform.

\begin{proposition}\label{propFredholm}
Let $P$ be a differential operator of weighted order $d$ on a closed foliated manifold $(M,H)$.
If the weighted principal symbol $\sigma_H(P)(x,\xi)$ is invertible,
then for every differential operator $A$ on $M$ of weighted order $\le d$,
there exists a constant $C>0$ such that
\[ \|Au\|\le C(\|Pu\|+\|u\|)\]
for every smooth function $u\in C^\infty(M)$.
The norms in the inequality are $L^2(M)$ norms.
\end{proposition}
{\bf Proof.}
Fix some point $m\in M$, and choose a neighborhood $U\subseteq M$ with foliation coordinates $x\in\RR^p\times\RR^q$ such that $x=0$ at $m\in U$. Having chosen these coordinates, we write
\[ P = P_m + \sum x_jQ_j+ R.\]
We explain the notation.
First, $P_m$ denotes the weighted principal part of $P$, with coefficients frozen at $m$, i.e.,
\[ P_m  = \sum_{||\alpha||=d} a_\alpha(m)\partial^\alpha.\]
Thus, $P_m$ is a constant coefficient operator on $\RR^n$, homogeneous for our weighted grading.
The remaining terms contain operators $Q_j$ of weighted order $d$, and a lower order part $R$.

Let $u$ be a smooth function in $U\subseteq M$ supported in a ball of small radius $|x|<\varepsilon$.
Then
\[ \|P_mu\|-\|Pu\|\le \|(P_m-P)u\|\le \varepsilon \sum \|Q_ju\| + \|Ru\|.\]
Taking Fourier transforms and applying the Plancherel theorem, one easily verifies that invertibility of $\sigma_H(P)$ implies the inequality
\[ \|Su\| \le C(\|P_mu\|+\|u\|),\]
for every constant coefficient operator $S$ of weighted order $\le d$ ($C$ depends on $S$, but not on $u$). Likewise
\[ \|Tu\| \le \varepsilon \|P_mu\|+C\|u\|,\]
for constant coefficient operators $T$ of weighted order $\le d-1$ (here $C$ depends on $T$ and $\varepsilon$).

Once these inequalities have been established for constant coefficient operators, the same inequalities for operators $S, T$ with smoothly varying coefficients follow directly (assuming $u$ is compactly supported). In particular, we have
\[ \|Au\|\le C(\|P_mu\|+\|u\|)\]
for $u$ supported in $|x|<\varepsilon$.
We also find,
\[ \|P_mu\|-\|Pu\|\le  \varepsilon \sum C_j(\|P_mu\| + \|u\|)+ \varepsilon \|P_mu\|+C_S\|u\|.\]
By taking $\varepsilon$ sufficiently small so that
\[ \varepsilon (1+\sum C_j)<\frac{1}{2}\]
we see that there exists a (large) $C>0$ such that,
\[ \|P_mu\|-\|Pu\|\le  \frac{1}{2}\|P_mu\| + C\|u\|\]
for $u$ compactly supported in the ball $|x|<\varepsilon$.
This, in turn, implies
\[ \|P_mu\|\le 2\|Pu\| + 2C\|u\|.\]
In summary, for every point $m\in M$ we can choose a neighborhood $V$ (corresponding to $|x|<\varepsilon$ for sufficiently small $\varepsilon$) such that the desired inequality
\[ \|Au\|\le C(\|Pu\|+\|u\|)\]
holds for all $u$ with support in $V$ (where $C$ now depends on $V$).
The global inequalities follow by a partition of unity argument and the compactness of $M$.

\QED

Proposition \ref{propFredholm} suggests the definition of a modified Sobolev norm.
Choose an open cover $\{U_j\}$ of $M$, and for each $U_j$ a set of vector fields $X^1_j,\ldots, X^n_j$ that are linearly independent at all points of $U_j$, and such that the first $p$ vector fields $X_j^1,\ldots, X_j^p$ span $H$. 
Let $\{\phi_j\}$ be a partition of unity subordinate to $\{U_j\}$.
Note that $\phi_j X_j^i$ is an operator of order 1 for $i=1,\ldots,p$, and of weighted order 2 otherwise.
For the positive integer $d$ we then define the {\em weighted Sobolev space} $W^d=W^d(M,H)$ on the foliated manifold $M$ as the completion of $C^\infty(M)$ with respect to the norm
\[ \|u\|^2_{W^d} = \sum_j\sum_{||\alpha||\le d}\,\|\phi_j X_j^\alpha u\|^2_{L^2(M)} \]
As usual, with this norm $W^d$ is a Hilbert space, and the equivalence class of the norm is independent of the choice of measure on $M$ or the choice of cover, partition of unity, or vector fields. 

This definition allows us to formulate the following corollary of Proposition \ref{propFredholm}.

\begin{corollary}\label{corClosure}
Let $P$ be a differential operator of weighted order $d$ on a closed foliated manifold $(M,H)$.
If the weighted principal symbol $\sigma_H(P)$ is invertible,
then the domain of the closure of $P$ is the weighted Sobolev space $W^d$,
and the a priori estimates of Proposition \ref{propFredholm} extend by continuity to all $u\in W^d$.
\end{corollary}
The Fredholm property of $P$ now follows from the a priori estimates. 

\begin{theorem}
Let $P$ be a differential operator of weighted order $d$ on a closed foliated manifold $(M,H)$.
If the weighted principal symbol $\sigma_H(P)$ is invertible,
then the closed operator $\bar{P}\colon L^2(E)\to L^2(F)$ is Fredholm.
\end{theorem}
{\bf Proof.}
The proof is a standard argument from elliptic theory.
We sketch the main steps here.
First, re-write the a priori estimates of $P$ as follows,
\begin{align*}
\|u\|_{W^d}^2 & \le C\,(\|Pu\|^2+\|u\|^2)\\
& = C\,\ang{(P^*P+1)u,u}\\
& = C\,\|(P^*P+1)^{1/2}u\|^2.
\end{align*}
In other words, the a priori estimates for $P$ are equivalent to the boundedness of the Hilbert space operator
\[ (P^*P+1)^{-1/2}\;\colon\; L^2(E)\to W^d(E).\] 
Observe that $W^d(E)$ is contained in the standard Sobolev space of order $k$, where $k$ is the largest integer such that $k\le d/2$. 
Since invertibility of $\sigma_H(P)$ implies that $P$ is at least of weighted order $2$,
we have $k\ge 1$.
It follows by the standard Rellich lemma that the operator $(P^*P+1)^{-1}$ is compact as an operator on $L^2(E)$. 

To prove that $P$ is Fredholm we must prove that $(\bar{P}P^*+1)^{-1}$ is also compact as an operator on $L^2(F)$. Since the weighted symbol of the {\em formal} adjoint $P^t$ is simply the matrix adjoint of $\sigma_H(P)(x,\xi)$, it is clear that $P^t$ has invertible symbol as well, and so by the same argument as above we deduce that $(P^{t*}P^t+1)^{-1}$ is compact on $L^2(F)$. This will finish the proof if we can show that the closure of the formal adjoint $P^t$ is identical to the Hilbert adjoint $P^*$.

To see that this is correct, the necessary analytical step is to establish that {\em on $L^2$-sections} one may identify the weak action of $P$ with the closure $\bar{P}$ of $P$. In other words, if $u\in L^2(E)$ and $v\in L^2(F)$ and $Pu = v$ weakly---in the sense of distribution theory---then it follows that $u\in W^d(E)$ and that $\bar{P}u = v$. (The proof is an application of Friedrichs mollifiers; see for example [Roe].) 
Since by definition $P^*u=v$ is equivalent to $P^tu=v$ weakly for $u, v\in L^2$, Corollary \ref{corClosure} implies that the closure of $P^t$ is indeed equal to the Hilbert space adjoint $P^*$ of $P$.

\QED

\section{The Tangent Groupoid for Foliations}

In this section we construct the tangent groupoid that is appropriate for the type of operators we are studying here. The construction is very similar to, but subtly different from the construction of the tangent groupoid for contact manifolds. For the original idea of the tangent groupoid method, see [Co].

As before, $(M,H)$ denotes a foliated closed manifold.
Algebraically, the tangent groupoid $\THM$ for $(M,H)$ is simply the union of smooth groupoids
\[ \THM = H\oplus N\;\cup\; M\times M\times (0,1].\]
Here $H\oplus N$ is conceived as a bundle of algebraically disjoint (graded) abelian groups,
while $M\times M\times (0,1]$ is a family of pair groupoids $M\times M$ parametrized by $t\in (0,1]$. 
As usual, the point is that the algebraic groupoid $\THM$ can be equiped with the structure of a smooth groupoid, by appropriately gluing together the two constituent pieces.

In order to achieve the desired gluing, we `blow up' the diagonal in $M\times M$ using the graded dilations
\[ \delta_t\;\colon\RR^{p+q}\to\RR^{p+q}\;\colon\;\delta_t(x,y)=(tx,t^2y)\]
in the first component of $M$.
To be more specific, for an open subset $U\subseteq M$, choose coordinates 
\[ (x,y)\;\colon\; U\to \RR^p\times \RR^q\]
compatible with the foliation, i.e., such that $\partial_1,\ldots,\partial_p$ span $H$.
Such coordinates induce coordinates $(x,y,\xi,\eta)\in \RR^p\times \RR^q\times \RR^p\times\RR^q$ on the total space of $H\oplus N$ (restricted to $U$).  
Then a local chart for $\THM$ is obtained by extending the coordinates $(x,y,\xi,\eta,0)$ on $H\oplus N$ to coordinates $(x+t^{-1}\xi,y+t^{-2}\eta,x,y,t)$ on $U\times U\times (0,1]$. 
One could say that the tangent groupoid makes rigorous the idea that $H\oplus N$ is an {\em infinitesimal} tubular neighborhood of the diagonal of $M\times M$,
with the added subtlety that we have modified the usual notion of the `order' of infinitesimals,
conform the weighted grading of our calculus.

It must be shown that different choices of local coordinates on $M$ induce a consistent smooth structure on $\THM$. This follows from a simple Taylor expansion.
Specifically, to study the effect of a change of coordinates on $M$, 
let $\phi\colon\RR^{p+q}\to\RR^{p+q}$ be a diffeomorphism that preserves the foliation structure and that fixes the origin, i.e., $\phi(0,0)=(0,0)$. 
The fact that $\phi$ fixes the foliation implies that there are smooth functions $f\colon\RR^{p+q}\to\RR^p$ and $g\colon\RR^q\to\RR^q$ such that
\[ \phi(x,y) = (f(x,y),\,g(y)).\]
The point is that $g$ is independent of $x$. A first order Taylor expansion gives
\begin{align*}
 \phi(tx,t^2y) &= (Df(tx,t^2y)+{\cal O}(t^2), Dg(t^2y)+{\cal O}(t^3))\\
 &= (tDf(x,0)+{\cal O}(t^2), t^2Dg(y)+{\cal O}(t^3)),
\end{align*}
which implies
\[ \delta_t^{-1}\phi\delta_t(x,y) = (Df(x,0),\,Dg(y)),\] 
which is the transformation law for coordinates on $H\oplus N$.
This simple fact guarantees that the smooth structure on $\THM$ is well-defined.
It also explains why we must take $H\oplus N$ as the groupoid at $t=0$, instead of $TM$.

\section{The Topological Index}

The tangent groupoid encodes, in a very nice way, the notion that the (noncommutative) algebra of operators on $M$ quantizes the (commutative) algebra of symbols on $T^*M$.

Recall that the $C^*$-algebra of a smooth groupoid $\GG$ is the completion 
of the convolution algebra $C_c^\infty(\GG)$ in a suitable norm.
The $C^*$-norm is defined as the supremum of the operator norms in the regular representations of $\GG$.
Following Connes' argument ([Co]), the restriction of smooth functions on $\THM$ to the $t=0$ fiber $H\oplus N$ induces a $\ast$-homomorphism,
\[ \pi_0\;\colon\; C^*(\THM)\to C^*(H\oplus N)\cong C_0(H^*\oplus N^*).\]
The kernel of this map is the contractible $C^*$-algebra
\[ C^*(M\times M\times (0,1]) \cong C_0((0,1])\otimes C^*(M\times M).\]
Therefore, the induced map in $K$-theory is an isomorphism
\[ \pi_0\;\colon\; K_0(C^*(\THM))\cong K^0(H^*\oplus N^*)\cong K^0(T^*M).\]
Restriction to the $t=1$ fiber induces the $\ast$-homomorphism
\[\pi_1\;\colon\; C^*(\THM)\to C^*(M\times M)\cong \KK(L^2(M)).\]
Combining the two maps, we obtain what we shall call the {\em topological index} for the foliation $(M,H)$,
\[ {\rm Ind}_H = \pi_1\circ\pi_0^{-1}\;\colon\; K^0(T^*M)\to K_0(\KK) = \ZZ.\]
Now, an invertible weighted principal symbol $\sigma_H(P)$
defines a compactly supported $K$-theory element
\[ [\pi^*E, \pi^*F, \sigma_H(P)] \in K^0(T^*M)\]
in exactly the same way that an elliptic symbol does.  
We have two vector bundles $\pi^*E$ and $\pi^*F$ over $T^*M$,
and an isomorphism $\sigma_H(P)$ between them that is defined outside a compact set.
Such a triple, by definition, determines an element in $K^0(T^*M)$.

We will prove two things.
First, we will show that our topological index ${\rm Ind}_H$ applied to the $K$-theory class $[\sigma_H(P)]$ computes the Fredholm index for the hypoelliptic operator $P$ (Theorem \ref{thm1} below).
Secondly, we will show that the topological index is in fact independent of the foliation (Theorem \ref{thm2}). 
Thus, the topological index ${\rm Ind}_H$ for hypoelliptic operators associated to a foliation is the same as the topological index computed by Atiyah and Singer for elliptic operators.
Theorem \ref{thm1} and Theorem \ref{thm2} together prove our index formula Theorem \ref{thmHormander}.

\section{The Index as a Graph Projection}

Before constructing the relevant elements in $K$-theory,
we quickly review the notion of a {\em graph projection}.
Let $T\in {\rm End}(V)$ be a linear endomorphism of a finite dimensional complex inner product space $V$.
The graph projection $e_T$ of $T$ is the orthogonal projection of vectors in $V\oplus V$ 
onto the graph $\{(v,Tv)\in V\oplus V\}$ of $T$. 
It is an elementary exercise in linear algebra to derive a formula for this projection as
\[
e_T = \left(\begin{array}{cc}(1+T^*T)^{-1}& (1+T^*T)^{-1}T^*\\T(1+T^*T)^{-1}& T(1+T^*T)^{-1}T^*\end{array}\right)\;
= \left(\begin{array}{cc}(1+T^*T)^{-1}& T^*(1+TT^*)^{-1}\\T(1+T^*T)^{-1}& 1-(1+TT^*)^{-1}\end{array}\right).
\]
More generally, if $T$ is a closed unbounded operator in a Hilbert space $\HH$,
then the graph projection of $T$ is a bounded operator on $\HH\oplus \HH$, 
and is computed by the same formula. 

To see the relevance to index theory of this construction,
consider a Fredholm differential operator $P\colon L^2(E)\to L^2(F)$
that satisfies a priori estimate as in Proposition \ref{propFredholm}.
The graph projection of $P$ is an operator on the Hilbert space $L^2(E)\oplus L^2(F)$.
The a priori estimates (together with the Rellich lemma in Sobolev theory) imply that
\begin{align*}
(1+P^*P)^{-1} &\in \KK(L^2(E)) \\
P(1+P^*P)^{-1}&\in \KK(L^2(E), L^2(F))\\
(1+PP^*)^{-1} &\in \KK(L^2(F))\\
P^*(1+PP^*)^{-1}&\in \KK(L^2(F), L^2(E)).
\end{align*}
Taken together these expressions amount to the statement that
\[ e_P - \left(\begin{array}{cc}0&0\\0&1\end{array}\right) \in \KK(L^2(E\oplus F)).\]
Thus, we can define a $K$-theory element
\[ [P] = [e_P] \ominus [\left(\begin{array}{cc}0&0\\0&1\end{array}\right)] \in K_0(\KK).\]
The significance of this construction is clear from the following proposition.
\begin{proposition}
Under the isomorphism 
\[ K_0(\KK)\cong \ZZ\]
that maps compact projections to their rank, 
the element $[P]$ constructed above corresponds to the Fredholm index of $P$.
\end{proposition}
{\bf Proof.}
The graph projections $e_t$ of the scaled operators $tP$, $t\ge 1$ form a norm continuous family, with $e_1 = e_P$.
As $t\to \infty$, this homotopy of projections converges in norm to the projection
\[ e_\infty = \left(\begin{array}{cc}[{\rm Ker}P]& 0\\0& 1-[{\rm Ker}P^*]\end{array}\right).\]
Here $[{\rm Ker}P]$ and $[{\rm Ker}P^*]$ denote the projections onto the kernels of (the closure of) $P$ and $P^*$, respectively. Thus, we have an equivalence of $K$-theory elements
\[ [e_P] \ominus [\left(\begin{array}{cc}0&0\\0&1\end{array}\right)] = [\left(\begin{array}{cc}[{\rm Ker}P]& 0\\0& 1-[{\rm Ker}P^*]\end{array}\right)] \ominus [\left(\begin{array}{cc}0&0\\0&1\end{array}\right)] =
[{\rm Ker}P] - [{\rm Ker}P^*].\]

\QED

\noindent The beauty of the preceding construction is that we can also apply it to the weighted principal symbol $\sigma_H(P)(x,\xi)$ of $P$.
Let $e_{\sigma_H(P)}$ denote the family of graph projections of $\sigma_H(P)$,
conceived as as a smooth projection-valued section in the bundle of $\ast$-algebras
\[{\rm End}(\pi^*E \oplus \pi^*F),\]
over $H^*\oplus N^*$.
We obtain an analytic $K$-theory element
\[ [\sigma_H(P)]=[e_{\sigma_H(P)}] \ominus [\left(\begin{array}{cc}0&0\\0&1\end{array}\right)] \in K_0(C_0(H^*\oplus N^*)).\]
The crucial point this time is that invertibility of $\sigma_H(P)$ for $\xi\ne 0$
guarantees that all the sections
\[
(1+\sigma^*\sigma)^{-1},\;
\sigma(1+\sigma^*\sigma)^{-1},\;
(1+\sigma\sigma^*)^{-1},\;
\sigma^*(1+\sigma\sigma^*)^{-1}
\]
vanish as $\xi\to\infty$. In other words, all these expressions represent $C_0$-sections over the locally compact space $H^*\oplus N^*$, so that
\[ e_{\sigma_H(P)} - \left(\begin{array}{cc}0&0\\0&1\end{array}\right) \in C_0(H^*\oplus N^*, {\rm End}(\pi^*E\oplus \pi^*F)),\]
as required.

\begin{proposition}
Under the canonical isomorphism between analytic and topological $K$-theory 
\[ K_0(C_0(H^*\oplus N^*)) \cong K^0(H^*\oplus N^*)\]
the element $[\sigma_H(P)]$ constructed above corresponds to the triple
\[ [\pi^*E, \pi^*F, \sigma_H(P)] \]
in compactly supported $K$-theory.
\end{proposition}
{\bf Proof.}
Let $G$ denote the range of the projection $e_{\sigma_H(P)}$ as a sub-bundle of $\pi^*E\oplus \pi^*F$,
and let 
\[ \tau\;\colon\; G\to \pi^*F\]
denote the restriction of the map $\pi^*E\oplus \pi^*F\to \pi^*F$ to $G$.
Clearly, $\tau$ is an isomorphism for $\xi\ne 0$.
 
The invertibility of $\sigma_H(P)$ implies that as $\xi\to\infty$ the fiber of $G$ converges
to the fiber of $\pi^*F$, and therefore $\tau$ converges to the identity. 
To see this, simply let $\xi\to\infty$ in the formula for the graph projection of $\sigma_H(P)$. 
Thus $G$ and $\pi^*F$ are {\em identical} sub-bundles of $\pi^*E\oplus \pi^*F$ at an appropriately chosen boundary of $H^*\oplus N^*$.
It follows that the formal difference 
\[ [e_{\sigma_H(P)}]\ominus[\left(\begin{array}{cc}0&0\\0&1\end{array}\right)]\]
in analytic $K$-theory
corresponds in compactly supported $K$-theory to the triple 
\[ [G,\pi^*F,\tau]\in K^0(H^*\oplus N^*).\]
Now $G$ is isomorphic to $\pi^*E$ by the restriction to $G$
of the projection $p\colon E\oplus F\to E$. 
Then $\sigma_H(P)=\tau\circ p^{-1}$ shows that the cycles $[G,\pi^*F,\tau]$ and $[(\pi^*E, \pi^*F, \sigma_H(P)]$ are isomorphic.
   
\QED

We see that essentially the {\em same} construction gives the Fredholm index of $P$ (as a $K$-theory class), and an element in topological $K$-theory $K^0(T^*M)$ for the principal symbol $\sigma_H(P)$.
Our next step is to further unite these two objects by considering the graph projection of a {\em single} operator $\PP$ on the tangent groupoid $\THM$ that encompasses both the operator $P$ and its principal symbol $\sigma_H(P)$.

\section{$K$-theoretic Proof of the Index Theorem}

The tangent groupoid provides the precise geometric context in which to combine the differential operator $P$ and its principal symbol $\sigma_H(P)$ into a single operator.
From this unified operator $\PP$ we can construct, by the same analytic method employed in the previous section, a single element in $K_0(C^*(\THM))$.
This single $K$-theory element is really a continuous deformation of the symbol class $[\sigma_H(P)]\in K^0(T^*M)$ to the Fredholm index of $P$ (as a class in $K_0(\KK)$).
The index theorem follows as an easy corollary.

In order to construct a differential operator $\PP$ on $\THM$ that smoothly connects the principal symbol (at $t=0$) with the operator $P$ (at $t=1$),
the principal symbol is best thought of as a smooth family $\{P_m\}$ of constant coefficient operators 
\[ P_m=\sum_{||\alpha||=d}\, a_\alpha(m)\partial^\alpha\]
on the fibers of the vector bundle $H_m\oplus N_m$.
In the language of smooth groupoids, the principal symbol corresponds 
to a {\em right-invariant family} on the groupoid $H\oplus N$.

If $\GG$ is a smooth groupoid with base $\GG^{(0)}$ and source and range maps $r, s\colon \GG\to\GG^{(0)}$, then a right-invariant family $T$ of differential operators is, by definition, a collection $T=\{T_b\}$ parametrized by base elements $b\in \GG^{(0)}$,
such that (1) each $T_b$ is a differential operator on the source fiber $\GG_b = s^{-1}(b)$;
(2) the coefficients of the family $T_b$ are smooth functions of $\GG$ (equivalently, $T$ is a differential operator on $\GG$)
and (3) the family is invariant under right-multiplication $R_\gamma\colon \GG_{r(\gamma)}\to \GG_{s(\gamma)}$ with elements $\gamma\in \GG$.

Applied to the groupoid $\GG = H\oplus N$ with base $\GG^{(0)}=M$, the fibers $\GG_m = s^{-1}(m)$ are just the fibers $H_m\oplus N_m$ of the vector bundle.
Then a right-invariant family is simply a collection of operators $T_m$, $m\in M$, 
such that each $T_m$ is a constant coefficient operator on $H_m\oplus N_m$,
with coefficients that are smooth functions on $M$.
Clearly, the principal symbol corresponds to such a right-invariant family.

The operator $P$ itself can also be conceived as a right-invariant family, 
this time on the pair groupoid $\GG=M\times M$.
The base of the pair groupoid is again $\GG^{(0)}=M$,
and each fiber $\GG_m = s^{-1}(m)=M\times \{m\}$ is simply a copy of $M$.
(We adopt the usual convention that arrows in the pair groupoid point from right to left,
so that $s(m,m')=m'$, $r(m,m')=m$.)
In this case, right invariance simply means that the operator is {\em the same} on each copy of $M$.
In other words, a right invariant family on $M\times M$ corresponds simply to a differential operator on the first factor $M$ in the Cartesian product.

Now the tangent groupoid is constructed precisely in such a manner that $P$ and $\{P_m\}$ can be combined into a single right-invariant family $\PP$ on the tangent groupoid $\THM$.
For elements $(m,t)\in M\times [0,1]$ in the base of $\THM$ we take 
\[ \PP_{(m,t)} = t^dP,\]
if $t\ne 0$, while for $t=0$ we let
\[ \PP_{(m,0)} = P_m.\]
The smooth structure on $\THM$ is constructed precisely in such a manner that the coefficients of this family $\PP$ are smooth functions on $\THM$.

The construction outlined so far is valid for any differential operator $P$
(and can easily be extended to the appropriate class of pseudodifferential operators).
The following proposition is the crucial ingredient in the proof of the index theorem,
and it relies on the invertibility of $\sigma_H(P)$ and the Fredholmness of $P$.

\begin{proposition}\label{propCont}
There exists an analytic $K$-cycle
\[ [\PP] \in K_0(C^*(\THM))\]
that restricts, at $t=0$, to the cycle
\[ [\sigma_H(P)]=[e_{\sigma_H(P)}] - [\left(\begin{array}{cc}0&0\\0&1\end{array}\right)]\in K^0(H^*\oplus N^*),\]
and that at $t=1$ restricts to
\[ [P]=[e_P] - [\left(\begin{array}{cc}0&0\\0&1\end{array}\right)]\in K_0(\KK).\]
\end{proposition}
The cycle $[\PP]$ is, of course, constructed from the graph projections $e_{(m,t)}$ for the operators $\PP_{(m,t)}$, by the method explained in the previous section. 
In other words, if we let 
\[ e_\PP = \{e_{(m,t)}, (m,t)\in M\times [0,1]\}\]
denote the entire family of graph projections, then 
\[ [\PP] = [e_\PP] \ominus [\left(\begin{array}{cc}0&0\\0&1\end{array}\right)].\]
The idea is that the family of projections
$e_{(m,t)}$ in the Hilbert spaces $L^2(\GG_{(m,t)})$
actually corresponds to the regular represention of a single element in (a matrix algebra over the unitalization of) the $C^*$-algebra $C^*(\THM)$.
It is remarkably tedious to prove this. 
The fact that $\PP$ is a smooth family is not sufficient.
It is, in particular, the continuity at $t=0$ of the family $e_{(m,t)}$ that is difficult to verify.
In some sense, this is the `hard nut' at the heart of the index theorem that is not yet cracked by the machinery presented here. 
In the last section of $[vE1]$ a detailed proof is presented, and we refer the interested reader to that paper.

An immediate corollary of Proposition \ref{propCont} is the $K$-theoretic version of the index theorem.

\begin{theorem}\label{thm1}
Let $P$ be a differential operator on a closed foliated manifold $(M,H)$
with invertible principal symbol $\sigma_H(P)$.
Then
\[ {\rm Index}\,P = {\rm Ind}_H([\sigma_H(P)]).\]
\end{theorem}
{\bf Proof.}
By definition, ${\rm Ind}_H$ is the composition of two maps: (1) the inverse of the isomorphism 
\[ K_0(C^*(\THM)\cong K_0(C^*(H\oplus N))\]
induced by restriction to the $t=0$ fiber in $\THM$, and (2) restriction at $t=1$.
But Proposition \ref{propCont} states precisely that $[\sigma_H(P)$ is the restriction at $t=0$ of a $K$-theory element in $K_0(C^*(\THM))$ that, in turn, restrict to the Fredholm index of $P$ at $t=1$.

\QED

All that is left  to prove is that the topological index ${\rm Ind}_H$ constructed here is identical to the topological index of Atiyah and Singer.
That proof is the content of the next and final section.

\section{The Computation of the Topological Index}\label{section:last}

We defined a topological index for a foliation $(M,H)$
\[ {\rm Ind}_H\;\colon\; K^0(H^*\oplus N^*)\to \ZZ\]
by means of the tangent groupoid $\THM$.
Any choice of section $N\to TM$ induces the same canonical isomorphism
\[ K^0(H^*\oplus N^*)\cong K^0(T^*M).\]
It is therefore hardly surprising that there is only one topological index.

\begin{theorem}\label{thm2}
The topological index 
\[ {\rm Ind}_H\;\colon\; K^0(T^*M)\to \ZZ\]
is independent of the foliation $H$.
\end{theorem}
{\bf Proof.}
We enlarge the parabolic tangent groupoid by introduction of a second parameter $s\in [0,1]$.
This larger groupoid is, in fact, the {\em adiabatic groupoid} of $\THM$.
For a general smooth groupoid $\GG$, the adiabatic groupoid is a groupoid fibered over $s\in [0,1]$ that `blows up' a tubular neighborhood of the space of units $\GG^{(0)}$ in $\GG$,
generalizing the way that Connes' tangent groupoid blows up the diagonal in $\GG=M\times M$.
(For a general discussion of this construction, see for example [Ni]). 

Recall that we can think of $\THM$ as a family of groupoids over the unit interval $[0,1]$, where at $t=0$ we have $H\oplus N$, while at each $t>0$ we have a copy of $M\times M$.
Algebraically, the adiabatic groupoid $\THM^{ad}$ is the union of a family of groupoids $\GG_{(t,s)}$ parametrized by $(t,s)\in [0,1]^2$, and defined as follows:
\begin{align*}
\GG_{(t,s)} &= M\times M,\; {\rm for}\;t>0,s>0,\\
\GG_{(t,0)} &= TM,\; {\rm for}\;t>0\\
\GG_{(0,s)} &= H\oplus N ,\; {\rm for}\;s\in [0,1].
\end{align*}
Since each groupoid $\GG_{(t,s)}$ has unit space $M$, the unit space of the adiabatic groupoid $\THM^{ad}$ is the manifold with corners $M\times [0,1]\times [0,1]$.  
Schematically:
\[ \xymatrix{{(t,s)=(0,0)}\ar@{~>}[r] & {H\oplus N} \ar@{.}[r]|-{\delta_t^{-1}}  \ar@{.}[d]_{s^{-1}} &  {TM} \ar@{.}[d]_{s^{-1}}^{\TM} \\
              & {H\oplus N}      \ar@{.}[r]^-{\THM}_-{\delta_t^{-1}} 
              &   {M\times M} & {(t,s)=(1,1)}\ar@{~>}[l]    }
\]
For a constant value of $s$, the `blow-up' along the $t$-axis is conform the graded dilations $\delta_t^{-1}$,
while for a constant value of $t$, the `blow-up' (following the general construction of an adiabatic groupoid) is simply by the factor $s^{-1}$. 
Observe that the $t=1$ edge contains a copy of the usual tangent groupoid $\TM$ of Connes.
In the present context it can be conceived as the tangent groupoid for the degenerate foliation $H=TM$.

The point of introducing the groupoid $\THM^{ad}$ is that it gives rise to a commutative diagram in $K$-theory, induced by restriction of functions on $\THM^{ad}$ to each of the four corners of the square $[0,1]^2$. We proceed step-by-step.

Restriction of elements in $C^*(\THM^{ad})$ to the $(t,s)=(0,0)$ corner,
\[ C^*(\THM^{ad}) \to C^*(\GG_{(0,0)})\cong C_0(H^*\oplus N^*), \]
induces an isomorphism in K-theory,
\[ K_0(C^*(\THM^{ad})\cong K^0(H^*\oplus N^*).\]
To see this, let $\GG_0$ denote the groupoid that is the union of the $t=0$ and $s=0$ edges in $\THM^{ad}$.
The restriction map $C^*(\THM^{ad})\to C^*(\GG_0)$ induces an isomorphism in $K$-theory,
because the kernel of this map is the contractible ideal $C_0((0,1]^2,\KK)$.
But $C^*(\GG_0)$ itself contracts to $C^*(\GG_{(0,0)})=C^*(H\oplus N)$. 

Now let $\alpha$ denote restriction to the edge $s=1$,
\[ \alpha\;\colon\; C^*(\THM^{ad}) \to C^*(\THM) ,\]
and $\beta$  restriction to the edge $t=1$,
\[ \beta\;\colon\; C^*(\THM^{ad}) \to C^*(\TM) .\]
Further restriction to the corner $(t,s)=(1,1)$ gives two $\ast$-homomorphisms,
\begin{align*}
\phi&\;\colon\; C^*(\THM) \to C^*(M\times M) \\
\psi&\;\colon\; C^*(\TM) \to C^*(M\times M).
\end{align*}
We obtain a commutative diagram,
\[ \xymatrix{  K_0(C^*(\THM^{ad})) \ar[d]^{\alpha} \ar[r]^{\beta} & K_0(C^*(\TM)) \ar[d]^{\psi} \\
               K_0(C^*(\THM)) \ar[r]_-{\phi}& K_0(C^*(M\times M)) }
\]
But $\phi$ is just our topological index ${\rm Ind}_H$,
while $\psi$ is the topological index for the degenerate case $H=TM$.
It is, therefore, just the topological index of Atiyah-Singer, and we denote it by ${\rm Ind}_{TM}$.
Moreover, a simple inspection of the definition of $\THM^{ad}$ shows that the maps $\alpha$ and $\beta$ induce the obvious isomorphisms in $K$-theory.

Thus, our diagram simplifies to,
\[ \xymatrix{  K^0(H^*\oplus N^*) \ar[d]^{\cong} \ar[r]^-{\cong} & K^0(T^*M) \ar[d]^{\rm Ind_{TM}} \\
               K^0(H^*\oplus N^*) \ar[r]_-{\rm Ind_H} & \ZZ }
\]
which shows that under the canonical isomorphism $K^0(H^*\oplus N^*)\cong K^0(T^*M)$
the topological index ${\rm Ind}_H$ is equal to the topological index ${\rm Ind}_{TM}$ for elliptic operators.

\QED

\noindent {\bf Remark.}
A full proof of Theorem \ref{thmHormander} of course requires that one computes the cohomological formula for the index map ${\rm Ind}_{TM}$. 
All we have shown here is that the formula for our class of hypoelliptic operators is the {\em same} as that for elliptic operators.
We can simply point to the third paper on index theory of elliptic operators by Atiyah and Singer for a computation of the elliptic formula ([AS3]).
Alternatively, an independent proof of this formula following the tangent groupoid methodology has been developed in [Hi].

\section*{References}

\def\item{\vskip2.75pt
plus1.375pt minus.6875pt\noindent\hangindent1em} \hbadness2500
\tolerance 2500
\markboth{References}{References}

\item{[AS1]}
M.\ F.\ Atiyah and I.\ M.\ Singer,
{\sl The index of elliptic operators I},
Ann.\ of Math.\ 87 (1968), 484–--530

\item{[AS3]}
M.\ F.\ Atiyah and I.\ M.\ Singer,
{\sl The index of elliptic operators III},
Ann.\ of Math.\ 87 (1968), 546–--604

\item{[BG]} R.\ Beals and P.\ Greiner,
{\sl Calculus on Heisenberg Manifolds},
Ann.\ of Math.\ Studies (119), Princeton, 1988.

\item{[Bo]} L.\ Boutet de Monvel, 
{\sl On the index of Toeplitz operators of several complex variables},
Invent. Math. 50 (1979), 249--272.

\item{[Co]} A.\ Connes,
{\sl Noncommutative Geometry},
Academic Press, 1994.



\item{[Hi]} N.\ Higson, 
{\sl On the $K$-theory proof of the index theorem},
Index theory and operator algebras (Boulder, CO, 1991), 67---86, 
Contemp.\ Math., 148, 1993.

\item{[H\"o]} L.\ H\"ormander,
{\sl On the index of pseudodifferential operators},
Elliptische Differentialgleichungen, Band II, 127–-146. Schriftenreihe Inst.\ Math.\ Deutsch.\
Akad.\ Wissensch.\ Berlin, Reihe A, Heft 8, Akademie-Verlag, Berlin, 1971.

\item{[Ni]} V.\ Nistor,
{\sl An index theorem for gauge-invariant families: the case of solvable groups},
Acta Math.\ Hungar.\ 99 (1--2) (2003), 155--183.


\item{[Roe]} J.\ Roe,
{\sl Elliptic operators, topology and asymptotic methods},
Second edition,
Pitman Research Notes in Math.\, 395.

\item{[Ta]} M.\ E.\ Taylor,
{\sl Noncommutative microlocal analysis, part I},
Mem.\ Amer.\ Math.\ Soc.\, vol. 313, AMS, 1984.

\item{[vE1]} E.\ van Erp, 
{\sl The Atiyah-Singer Index Formula for Subelliptic Operators on Contact Manifolds, Part I},
Ann.\ of Math., to appear

\item{[vE2]} E.\ van Erp, 
{\sl The Atiyah-Singer Index Formula for Subelliptic Operators on Contact Manifolds, Part II},
Ann.\ of Math., to appear

\end{document}